\documentclass[11pt,a4paper]{article}

\usepackage[utf8]{inputenc}
\usepackage{amsmath}
\usepackage{amsfonts}
\usepackage{amssymb}
\usepackage{makeidx}
\usepackage{graphicx}
\graphicspath{ {./} }
\usepackage{listings}
\usepackage{tocloft}
\usepackage{verbatim}
\usepackage{xcolor}

\usepackage{calc}
\usepackage{CJKutf8}
\usepackage[hidelinks]{hyperref}

\newtheorem{le[mm]a}[theorem]{Le[mm]a}

\usepackage{arabtex}

\newcommand{\qed}{\nobreak \ifvmode \relax \else
      \ifdim\lastskip<1.5em \hskip-\lastskip
      \hskip1.5em plus0em minus0.5em \fi \nobreak
      \vrule height0.75em width0.5em depth0.25em\fi}
      
       \usepackage[polutonikogreek,english]{babel}
    
\usepackage[normalem]{ulem}

\usepackage{tikz, tkz-tab, pgfkeys}
\usetikzlibrary{calc,decorations.pathmorphing,shapes}

  \newcounter{exer}[section]

    \makeindex  
     
\title{Spherical trigonometry before the modern era: \\
  The treatise of Na\d{s}\=\i r al-D\=\i n al-\d{T}\=us\=\i }
\author{Athanase Papadopoulos}
 \date{\today}

\begin{document}

\maketitle

 \bigskip

\noindent {\bf Abstract} 
 
This is an overview of Na\d{s}\=\i r al-D\=\i n al-\d{T}\=us\=\i 's  \emph{Treatise of the quadrilateral}, an invaluable 13th century document on spherical geometry which was translated into French in 1891. The title we are using here is the one given by the translator (Alexandre Carath\'eodory). A title which is closer to the original Arabic is \emph{Disclosing the secrets of the secant figure}. The term   ``secant figure", to which the title refers, is the so-called ``complete (spherical) quadrilateral", that is, the figure that underlies what we call today Menelaus' Theorem. This theorem gives a formula that was extensively used by astronomers in their computations and the establishment of their tables since the first century AD, notably by Ptolemy, in the absence of the spherical trigonometric formulae that were discovered later. Na\d{s}\=\i r al-D\=\i n al-\d{T}\=us\=\i 's  treatise contains much more than Menelaus' theorem, since we find there a complete system of spherical trigonometric formulae, with complete proofs. The treatise includes at the same time invaluable historical information on the discovery of the trigonometric formulae by the Arab mathematicians of the Middle-Ages and the transformation of the field of spherical trigonometry that this discovery led to.

The final version of this paper will appear in the book \emph{Spherical geometry in the eighteenth century, I: 
 Euler, Lagrange and Lambert}, edited by
 Renzo Caddeo and Athanase Papadopoulos, Springer, 2026.

            

  \section{Introduction}

Na\d{s}\=\i r\index{Tusi@al-\d{T}\=us\=\i , Na\d{s}\=\i r al-D\=\i n} al-D\=\i n al-\d{T}\=us\=\i 's\footnote{Na\d{s}\=\i r al-D\=\i n al-\d{T}\=us\=\i \ was  a mathematician, astronomer and  philosopher. He belongs to the last period of the Arabic golden era in science. He was born in 1201 in the city of
\d{T}\=us in the province of Khoras\=an (North-East of Iran). We owe to him commented Arabic editions of several important Greek mathematical texts, including geometrical works of\index{Autolycus of Pitane} Autolycus of Pitane, Euclid,\index{Euclid} Apollonius,\index{Apollonius of Perga} Archimedes,\index{Archimedes} Theodosius,\index{Theodosius of Tripoli} Menelaus\index{Menelaus of Alexandria} and Ptolemy\index{Ptolemy of Alexandria}. These editions played an important role in the transmission of the Greek geometric works to Europe. An exposition of Euclid's \emph{Elements} attributed to him (the \emph{Book of exposition of the \emph{Elements}}), published in Rome in 1594, contains a commentary that was used and quoted by John Wallis and by Girolamo Saccheri in their works on the parallel problem, see \cite{R1, Rosenfeld, RY}. 
 Al-\d{T}\=us\=\i \  also wrote a treatise on the theory of parallels, called \emph{Al-ris\=ala al-sh\=afiya `an al-shakk f\=\i \ l-khu\d{t}\=u\d{t}} (the letter that heals the doubt raised by the parallel lines). The treatise is analyzed by Rosenfeld in \cite[p. 74-80]{Rosenfeld}. We also mention al \d{T}\=us\=\i's contribution on the theory of parallels in Chapter 7 of the present volume, written by the author of the present chapter together with Guillaume Th\'eret \cite{Papa-Theret-Lambert}.
 Al-\d{T}\=us\=\i\  also wrote treatises on logic, astronomy, algebra, number theory and combinatorics, see \cite{Carra, K1, K2, Ragep, RR}.} \emph{Treatise of the quadrilateral}\footnote{We are using the name that the translator, Alexandre Carath\'eodory (1833-1906), gave to this treatise, in French, \emph{Trait\'e du quadrilat\`ere} \cite{Tusi}. 
 The complete title, in Arabic, is \index{Carath\'eodory, Alexandre}  
 \emph{Kashf al-qin\=a` `an asr\=ar al-shakl al-qa\d{t}\d{t}\=a`}  (Disclosing the secrets of the secant figure). The work became known in Europe as \emph{Trait\'e du
quadrilat\`ere}  (Treatise of the quadrilateral), after Carath\'eodory's translation. Carath\'eodory 
  belonged to a distinguished Greek family from Constantinople, who was close to the Sublime Porte. His father, Stephanos  Carath\'eodory,\index{Carath\'eodory, Stephanos} was the personal physician of the Sultan Mahmud II. Alexandre Carath\'eodory studied law in Paris and obtained there a doctorate, with a dissertation on the philosophical side of  theory of error and its applications in the sciences. After his return to Constantinople, he  became an Ottoman civil servant and started a diplomatic career. 
He became ambassador of the Ottoman empire to Rome and the main negotiator at the Preliminary Peace Treaty of San Stefano which ended the
 Russo-Turkish War (1877-78). He was the great-uncle of the mathematician\index{Carath\'eodory, Constantin} Constantin Carath\'eodory.} is an invaluable 13th century treatise on spherical geometry. It was edited  and translated into French by Alexandre Carath\'eodory in 1891 \cite{Tusi}.
  The ``quadrilateral" to which the French title refers is a complete (Euclidean or spherical) quadrilateral; see Figures \ref{fig:EuclideanSector} and \ref{fig:Quadrilaterals-Tusi1}. It also goes by the names ``sector figure", ``secant figure"\index{secant figure} and\index{sector figure} others. We shall explain these terms below.
 
 Al-\d{T}\=us\=\i 's treatise is a testimony of the transformation of the subject of spherical geometry by the Arab mathematicians that started in the 9th century, when the use of the sine rule and the other formulae of spherical trigonometry replaced that of the sector figure. 
The exposition of a complete set of spherical trigonometric formulae is carried out in the last part of this work (Book V). The author establishes there the main formulae that are the subject of every modern treatise on spherical trigonometry, including the notion of polar triangle and its use in spherical trigonometry.    At the same time, Na\d{s}\=\i r al-D\=\i n's treatise constitutes a remarkable historical document that makes an overview of the state of the art of the field of spherical geometry as it was developed by the Arabs between the tenth and the thirteenth centuries.

 The treatise is divided into five books, with the following titles:

\begin{itemize}
\item {Book I}. On compounded ratios and their rules.\index{compounded ratio}

\item {Book II}. On the plane quadrilateral figure and the ratios that we find therein.

\item {Book III}. On the lemmas for the spherical quadrilateral figure and on what is necessary for
using it successfully.

\item {Book IV}. On the spherical quadrilateral figure and the ratios that we find
therein.

\item {Book V}. An exposition of the methods that dispense with the theory of the quadrilateral figure,  for what concerns the
knowledge of the arcs of great circles.

\end{itemize}

 The details of the content of this treatise are given in the rest of the present chapter.
 
   There is a review of a small part of  Na\d{s}\=\i r al-D\=\i n's treatise in Rosenfeld's  book on the history of non-Euclidean geometry \cite[p. 20ff]{Rosenfeld}. Rosenfeld writes there that this book contains   ``the most complete exposition of spherical geometry in the East in the Middle Ages, written by the greatest mathematician and astronomer of the 13th century."  \cite[p. 20]{Rosenfeld}

For the general mathematician, the title of Book I of the treatise deserves an explanation, since it is not clear a priori why  an object like a compounded ratio\index{compounded ratio} (a product of two ratios of numbers or magnitudes, that is, an expression of the form $\frac{A}{B}\cdot\frac{C}{D}$), needs\index{compounded ratio} a special chapter in such a book. We shall explain this in \S \ref{s:compounded}.  
  
  All along this chapter, we shall use modern notation, in particular for what concerns trigonometric functions, as Carath\'eodory does in his edition. 
 
 The first three sections that follow the present introduction are expositions of interesting remarks made by  al-\d{T}\=us\=\i \ in his treatise, on the intersection of three or four lines in the plane or on the sphere.
 
 In Section \ref{s:compounded}, we give some explanations on compounded ratios, which is the subject of the first book of  Al-\d{T}\=us\=\i 's treatise. 
 
Section \ref{s:sector-theorem} is concerned with the Sector Figure,  that is, Menelaus' Theorem,\index{Menelaus theorem} which\index{Theorem!Menelaus} is extensively used and thoroughly discussed in Al-\d{T}\=us\=\i 's treatise.

In Section \ref{s:summary-proposition}, we make a quick summary of Book V of the treatise, which is interesting in many aspects. We quote several series of important propositions from this book, and we linger on the importance of the treatise for the historical information it contains.

        \section{On the figures formed by the intersection of three and four lines in the Euclidean plane}\label{s:Figures}

\begin{figure}[htbp]
\centering
\includegraphics[width=.70\linewidth]{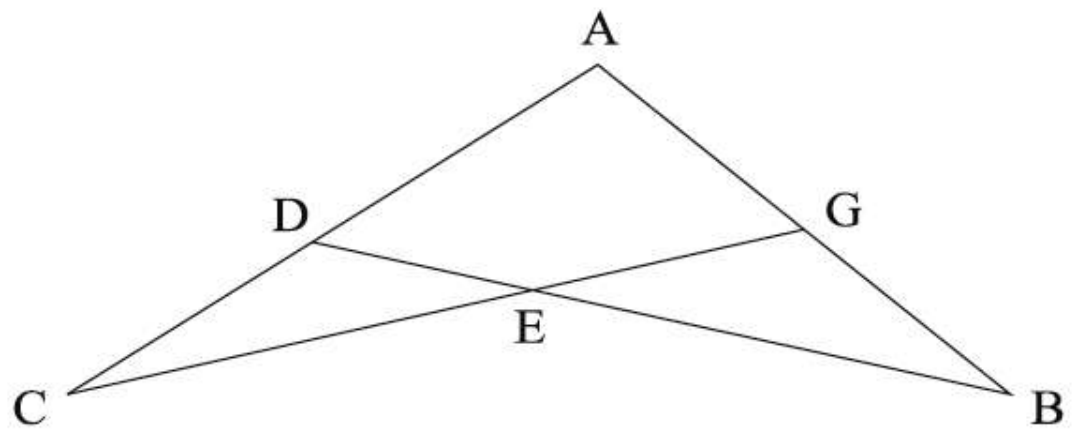}
\caption{\small {The Euclidean sector figure, or complete Euclidean quadrilateral}}
\label{fig:EuclideanSector}
\end{figure}

Among the simplest elements of Euclidean plane geometry, there are points, next, pairs of points, and the latter define lines. Then, one may consider a pair of two lines which, if they are generic (that is, if they are not
coincident or disjoint), define an angle. Going forward, the next step is to consider a 
triple of lines, which, again if they are generic, define a triangle. One may
continue in the same manner and make a further step, considering the union of
four generic lines. The resulting figure contains a subfigure represented in Figure \ref{fig:EuclideanSector}, made up of a quadrilateral to which are appended two triangles.
 This is the so-called \emph{(Euclidean) sector figure}, also called \emph{complete (Euclidean) quadrilateral}.

 The study of the sector figure arises naturally in Euclidean geometry. It comes naturally next to that of triangles. The\index{Theorem!Menelaus} classical Menelaus Theorem,\index{Menelaus theorem}  which we recall in \S \ref{s:spherical-quadr}, concerns the sector figure.

 Now we pass to spherical geometry: first, the intersection of three great circles and then, the intersection of four. 
 
  \section{On the figure formed by the intersection of three great circles}\label{s:three-circles-spherical}

 \begin{figure}[!ht]
\centering
 \includegraphics[width=0.5\linewidth]{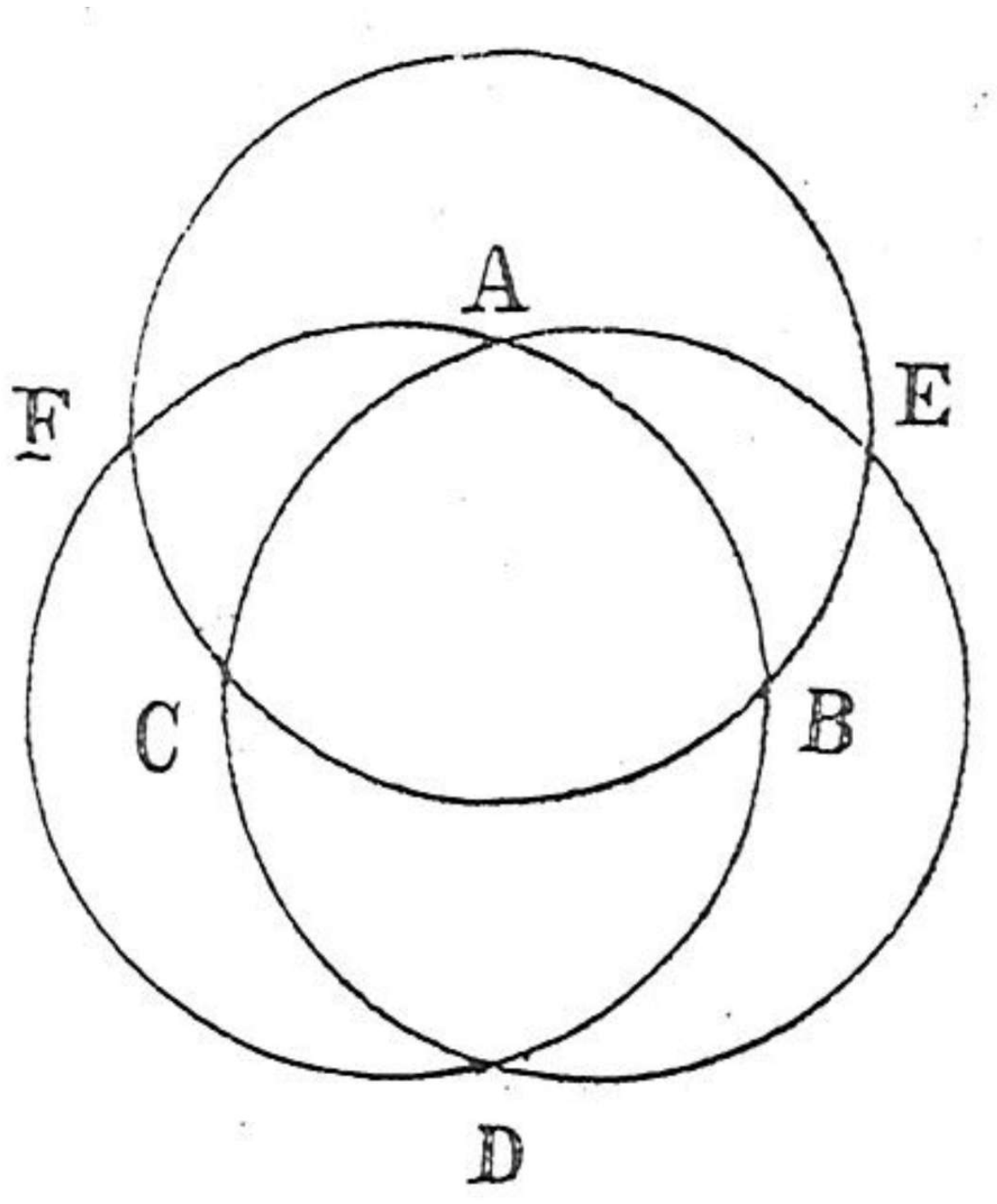}    \caption{\small 
The eight triangles made by three great circles on the sphere. The figure is extracted from Carath\'eodory's edition of al-\d{T}\=us\=\i 's treatise. Note that there is a triangle with sides $EF, FD, DE$.  It is remarkable that Na\d{s}\=\i r al-D\=\i n al-\d{T}\=us\=\i \ drew an abstract (non-realistic) figure of the intersection pattern of the three great circles of the sphere.}   \label{fig:Trois-cercles}  
\end{figure}

A spherical triangle is determined by three distinct points on the sphere together with the three arcs of great circles that join them pairwise.  Before studying spherical  triangles, Na\d{s}\=\i r al-D\=\i n   studies the intersection of three great circles (Chapter 2). He notes that a generic triple of great circles divides the sphere into eight spherical triangles (Figure \ref{fig:Trois-cercles}).  These eight triangles are pairwise equal,\footnote{Following Euclid, we use the term ``equal" to means the modern term ``isometric".} since they are pairwise images of each other by the antipodal map of the sphere. For instance, the triangles $ACF$ and $DEB$ of Figure \ref{fig:Trois-cercles} are equal. Thus, we have, in this figure,  at most four unequal triangles. Among these four triangles any two have at least one side and one angle in common, and the rest of their elements are supplementary of each other (their sum is equal to two right angles).\footnote{We recall that the length of a spherical segment is the angle made by the two rays that join the origin of the sphere to the two extremities of this segment.} For example, comparing the triangles $ABC$ and $CDF$, we see that the two angles at $C$ are equal, $AB=FD$, $BC$ and $FC$ are supplementary, $AC$ and $CD$ are supplementary, the angle $ABC$ is equal to the angle $AFC$ (they are the two angles made by two half great circles with common endpoints) and is supplementary to the angle $CFD$, and the angle $BAC$ is equal to the angle $BDC$ and is supplementary to the angle $CDF$. 

Now recall that in spherical geometry, triangles are determined by the measures of their angles.\footnote{This is Proposition 19 of Menelaus' \emph{Spherics}, see \cite[p. 548]{RR2}.} Thus, if we know the triangle $ABC$, then we know the triangle $CDF$.  The same reasoning applies to all the others triangles of the figure, pairwise.  
 Consequently, if one of the eight triangles in this figure is known, then all the other seven triangles are also known.

\section{Spherical complete quadrilaterals}\label{s:spherical-quadr}

We have already mentioned complete Euclidean quadrilaterals (see Figure \ref{fig:EuclideanSector}). In this section, we examine complete spherical quadrilaterals.
Thus, after spherical triangles, we study spherical quadrilaterals, and more especially \emph{complete quadrilateral}.\index{complete quadrilateral!spherical}
 
\begin{figure}[htbp]
\centering
\includegraphics[width=.60\linewidth]{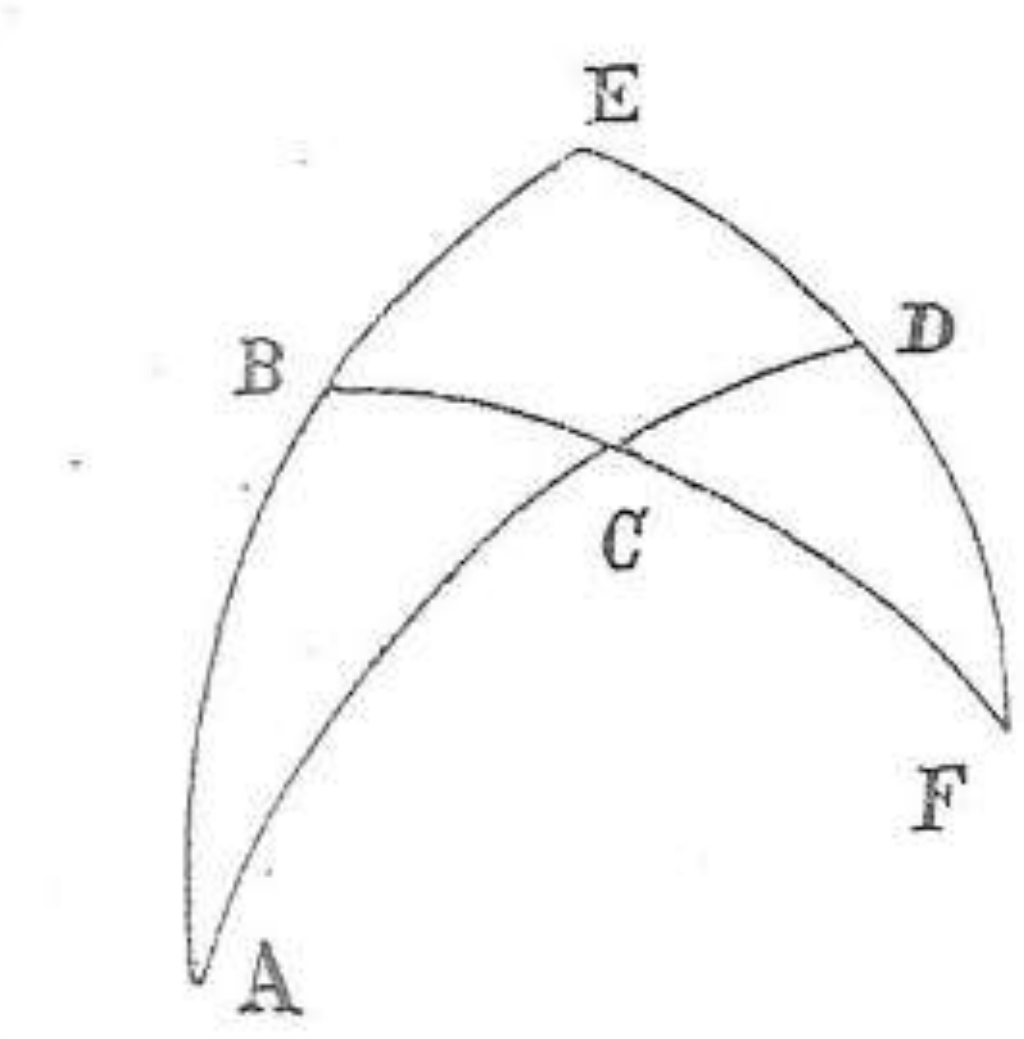}
\caption{\small A complete spherical quadrilateral. Figure extracted from Na\d{s}\=\i r al-D\=\i n al-\d{T}\=us\=\i 's treatise.}
\label{fig:Quadrilaterals-Tusi1}
\end{figure}

A \emph{spherical quadrilateral} is a convex figure on the sphere bordered by four arcs of great circles (the sides of the quadrilateral).
  A spherical quadrilateral becomes a \emph{complete spherical quadrilateral} if one adds to it
  two spherical triangles by extending two pairs of opposite sides until in  two sides in each such a pair meet.  
  Figure  \ref{fig:Quadrilaterals-Tusi1} is a complete spherical quadrilateral, extracted from Na\d{s}\=\i r al-D\=\i n al-\d{T}\=us\=\i 's treatise.

\begin{figure}[htbp]
\centering
\includegraphics[width=.60\linewidth]{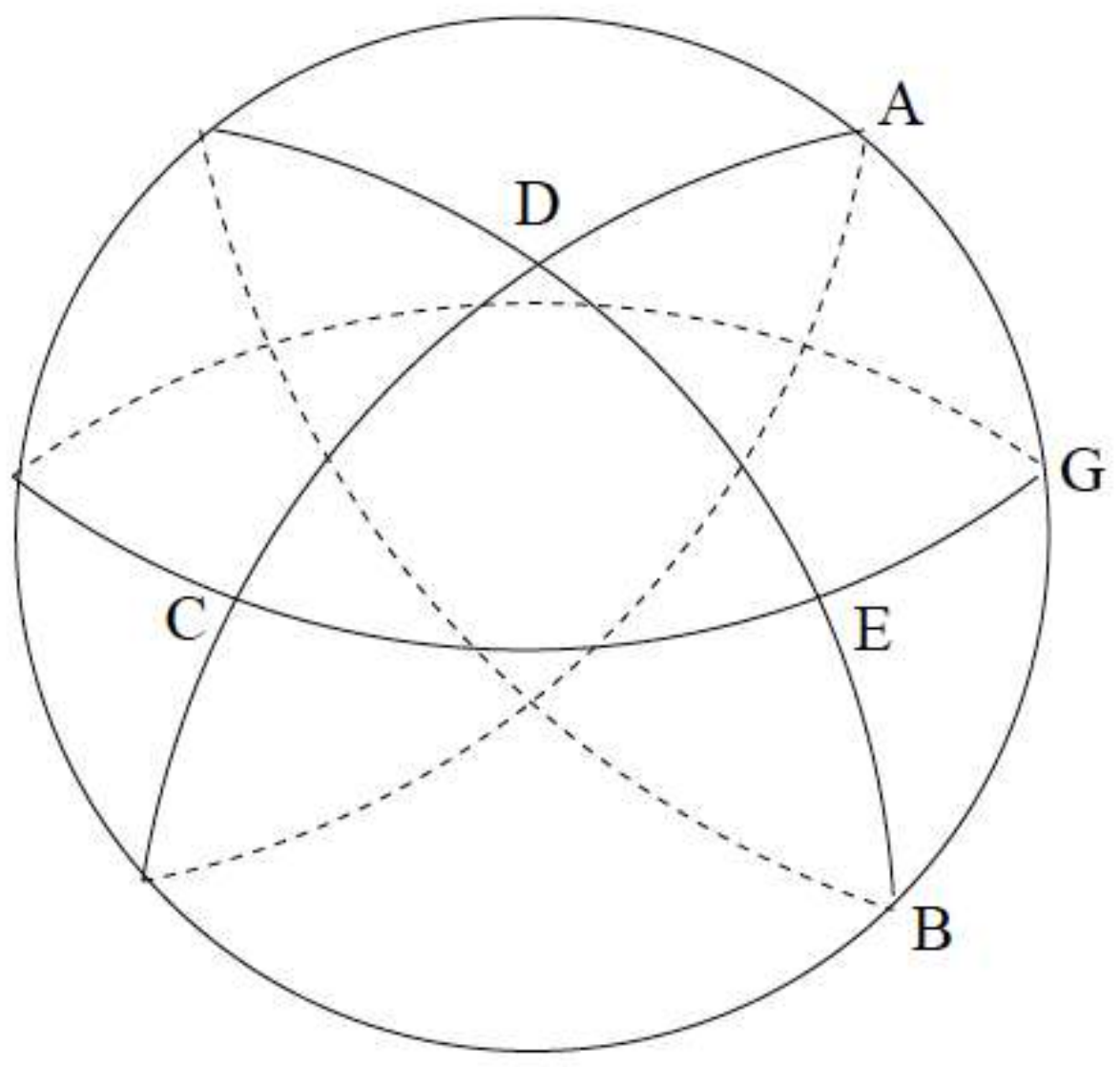}
\caption{\small The spherical quadrilateral $ADEG$ to which are appended the triangles $DCE$ and $EGB$ is a complete spherical quadrilateral.}
\label{fig:Quadrilaterals-Tusi}
\end{figure}

  In Figure \ref{fig:Quadrilaterals-Tusi} 
  we have have drawn four great circles, namely, the three great circles containing the sides of the spherical triangle $CDE$, and the fourth one, passing by the points $A$ and $G$.
  In this figure, also reproduced from a figure in Na\d{s}\=\i r al-D\=\i n al-\d{T}\=us\=\i 's 's treatise,
  the quadrilateral $ADEG$ to which are appended the triangles $DCE$ and $EGB$ is a complete spherical quadrilateral.  Note that each spherical  quadrilateral gives rise to four complete spherical quadrilaterals, depending on the choice of the two adjacent sides and the direction in which we extend them. In Figure \ref{fig:Quadrilaterals-Tusi}, made by the intersection of the three circles, we have a
total of six quadrilaterals which are equal (by the antipodal map), which makes at most 
 three non-equal quadrilaterals. Since each quadrilateral
gives rise to four sector figures, we have in all $6 \times 4 = 24$ sector figures. But we have at most 12 non-equal ones.

The study of complete spherical quadrilaterals\index{complete quadrilateral!spherical} was carried out extensively  by  Menelaus of Alexandria\index{Menelaus of Alexandria} in his \emph{Spherics}. 
   This figure is associated with a proposition of this treatise which had far-reaching applications and which is known in the modern mathematical literature (although stated in a slightly different form) as \emph{Menelaus theorem}, see \cite[Proposition 66]{RR2}. The\index{Menelaus theorem} proposition\index{Theorem!Menelaus}  was used by Ptolemy\index{Ptolemy of Alexandria} in Book I, Chapter 11 of the \emph{Almagest}, the famous astronomical treatise written around the middle of the second century AD, as the main lemma for the astronomical calculations that are carried out there, see \cite[p. 50--55]{Halma}. We have elaborated at length on this proposition and its applications in the commentary contained in our edition of the \emph{Spherics}, see \cite[p. 292--337]{RR2}.

 The complete spherical quadrilateral, together with the proposition to which it belongs, also bears the names ``quadrilateral figure", ``complete quadrilateral",\index{complete quadrilateral!spherical}``sector figure" or ``secant figure",\index{secant figure} and\index{sector figure} there are others. The name ``secant figure"   was used by Ptolemy in the \emph{Almagest}, and it is also the name that was adopted by the Arab mathematicians\footnote{The name \RL{Al^sakl Alqa.t.taa`}, meaning the ``secant figure", or the ``cutting figure", was used by Th\=abit ibn Qurra (826--921) and others.} in their commentaries on Menelaus' and\index{Menelaus theorem} Ptolemy's works\index{Theorem!Menelaus}  and in the treatises on spherical geometry they wrote. The title of 
  Na\d{s}\=\i r al-D\=\i n's treatise, which is the subject of the present chapter, refers to the complete spherical quadrilateral.

\section{Compounded ratios}\label{s:compounded}

Book I of  Na\d{s}\=\i r al-D\=\i n's treatise  is concerned with 
compounded ratios. In modern notation, a compounded ratio\index{compounded ratio} is
an expression of the form $\frac{A}{B}\cdot \frac{C}{D}$. 
Menelaus' proposition called the Sector Figure, which we shall discuss in the next section and the rest of the treatise, make heavy use of such expressions.  Therefore, a skill in the manipulation of compounded ratios is needed in following the arguments.
   
We remind the reader that expressions involving compounded ratios, such as the equality $\frac{A}{B} \cdot\frac{C}{D}= \frac{A}{D}\cdot \frac{C}{B}$, that are obvious for us, were not so for the Ancients. The reason has to do with the differences in our points of view on such expressions and the operations they involve. Indeed, we are used (consciously or unconsciously) to consider division of quantities and multiplication as formal algebraic operations pertaining to a certain mathematical structure (ring or field), and the identities that they satisfy follow from the axioms of this structure: commutativity, distributivity, etc. For a mathematician of Greek antiquity or of the Arabic period that followed, this was not the way mathematics was conceived. Their point of view goes back to Euclid's \emph{Elements}, where the operations of multiplication and division have a geometric meaning, and the manipulation of expressions involving compounded ratios and equalities between them needed proofs. 

%
%
 
As a matter of fact, there is hardly any definition of compounded ratios\index{compounded ratio} in Euclid's \emph{Elements}. The first allusion to them is in Definition 5 of Book VI, which is considered to be a later addition; see Heath's comments in \cite[Vol. 2, p. 189]{Euclid-Heath}. The definition, in Heath's translation, says the following: \emph{A ratio is said to be compounded of ratios when the sizes of the ratios multiplied together make some (? ratio or size).}\footnote{An uncertainty is expressed by the question mark in Heath's edition.}

Proposition 23 of Book VI of the \emph{Elements} is usually referred to as a typical place where compounded ratios are implicitly used. The proposition is purely geometric and it concerns areas of parallelograms. It says the following \cite[Vol. 2, p. 247]{Euclid-Heath}: \emph{Equiangular parallelograms have to one another the ratio compounded of the ratios of their sides}. It is interesting to read Heath's notes on this proposition.  In his \emph{History of Greek Mathematics}, he writes  \cite[Vo. 1, p. 393]{Heath-History3} (see Figure \ref{fig:Heath23} here):

 \begin{figure}[!ht]
\centering
 \includegraphics[width=0.6\linewidth]{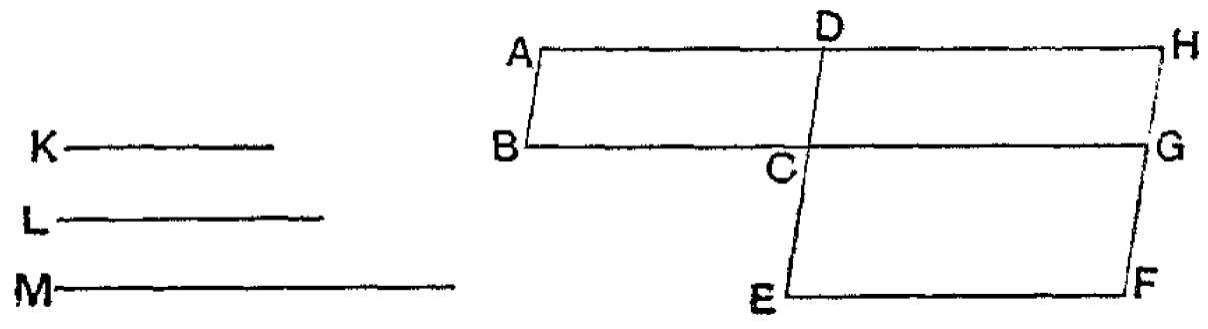}    \caption{\small Figure for Euclid's Proposition 23 of Book VI (from Heath's  \emph{History of Greek Mathematics})}   \label{fig:Heath23}  
\end{figure}
\begin{quote}\small

Proposition 23 is important in itself, and also because it introduces us to the practical use of the method of compounding, i.e. multiplying, ratios which is of such extraordinarily wide application in Greek geometry. Euclid has never defined ``compound ratio" or the ``compounding" of ratios; but the meaning of the terms and the way to compound ratios are made clear in this proposition. The equiangular parallelograms are placed so that two equal angles as $BCD, GCE$ are vertically opposite at $C$. Complete the parallelogram $DCGH$. Take any straight line $K$, and find another $L$, such that  \[BG:CG=K:L,\]
and again another straight line $M$, such that
\[DC:CE=L:M.\]
Now the ratio compounded of $K:L$ and $L:M$ is $K:M$; therefore $K:M$ is the ``ratio compounded of the ratios of the sides".
And 
\[(ABCD):(DCGH)=BC:CG=K:L;\]
\[(DCGH):(CEFG)=DC:CE=L:M,\]
therefore, \emph{ex aequali},
\[(ABCD):(CEFG):K:M.\]
\end{quote}

 This excursus should explain why Na\d{s}\=\i r al-D\=\i n makes a systematic\index{compounded ratio} study of compounded ratios in his treatise. In Book I, he gives 11 propositions on this subject. As an example, we quote one of them, namely,   Proposition 10. It says the following (in modern terms):
 
 \begin{quote}\small
 
 If we have an equality between a ratio of two quantities and a compounded ratio of four quantities expressed as $\frac{A}{B}=\frac{C}{E}\cdot\frac{D}{F}$, 
then the ratio of any one of the three quantities $A,E,F$ to any one of the three quantities $B,C,D$ is equal to a compounded ratio of the four remaining quantities, whose numerators are taken among the remaining quantities  in $A,E,F$ and the denominators among the remaining quantities in $B,C,D$.

\end{quote}

Thus, this proposition says in particular that 
 if $\frac{A}{B}=\frac{C}{E} \cdot \frac{D}{F}$, then $\frac{A}{C}=\frac{B}{E}  \cdot\frac{D}{F}$.  Na\d{s}\=\i r al-D\=\i n gives a proof of this proposition which is based on Propositions 33 and 34\footnote{We are using the numbering in Heath's edition. Na\d{s}\=\i r al-D\=\i n's numbering is different.} of Book XI of the \emph{Elements}, which say respectively \cite[Vol. 3, p. 342 and 345]{Euclid-Heath} that
\emph{Similar parallepipedal solids are to one another in the triplicate ratio of their corresponding sides} and
 \emph{In equal parallelepipedal solids the bases are reciprocally proportional to the heights; and those parallelepipedal solids in which the bases are reciprocally proportional to the heights are equal.} Thus, the language for compounded ratios is still geometric.

Now we pass to the rest of Na\d{s}\=\i r al-D\=\i n's treatise, where compounded ratios\index{compounded ratio} are extensively used.
  
 \section{The sector figure theorem}\label{s:sector-theorem}
 
We start by a few words of explanation on the expression ``sector
figure". This expression refers to a proposition (or a formula) and the corresponding 
figure. There is a Euclidean and a spherical sector figure.\footnote{In Arabic, the same word is used for Figure and Proposition, ``shakl" (\RL{^sakl}). This is not unrelated to the fact that in Euclid's \emph{Elements}, in Theodosius' \emph{Spherics},  in Menelaus' \emph{Spherics} and in other geometrical Greek treatises, to each proposition corresponds one and only one figure. Thus, the references to some proposition of one of these treatises is at the same time a reference to the corresponding figure.}

We start with the Euclidean proposition. The associated figure is represented in Figure \ref{fig:EuclideanSector}, in which we have a convex\index{secant figure} quadrilateral\index{sector figure} $AGED$ whose two pairs of
opposite sides $AG$ and $DE$ are produced until they intersect. In other words, to the quadrilateral $AGED$ in Figure \ref{fig:EuclideanSector}, we adjoin the two triangles
CED and BEG. 

%
%
%
With this notation, the Euclidean sector figure says that
$$\frac{AG}{BG}=
\frac{AC}{CD}\cdot \frac{DE}{EB}.
$$

The use and the proof of this proposition also include the cases where one or two opposite edges AD and GE are parallel.

There is a similar construction on the sphere which gives a \emph{spherical sector figure}.
  Using the same notation as in Figure \ref{fig:EuclideanSector} and assuming now that the quadrilateral is spherical, this proposition says:
$$\frac{\sin AG}{\sin BG}=
\frac{\sin AC}{\sin CD}\cdot \frac{\sin DE}{\sin EB}.
$$
%
%
%

%
%
%
This is Menelaus' Theorem on the sphere, proved in \cite{RR2} and used by Ptolemy in the \emph{Almagest}. It is cumbersome, compared to the sine rule. This is why the Arabs, when they discovered the sine rule, called it sometimes ``the proposition which exempts" (that is, it exempts from using Menelaus' Theorem). Na\d{s}\=\i r al-D\=\i n al-\d{T}\=us\=\i \ uses this name as well.

 \section{A summary of Book V and some propositions}\label{s:summary-proposition}
  Book V of Na\d{s}\=\i r al-D\=\i n's treatise contains seven chapters which are concerned with spherical trigonometry.
The author studies in particular with the following  topics:

\begin{itemize}
\item (Chapter 2) A complete classification of the figures that are formed by three great circles on the sphere. 
\item (Chapter 3) A complete classification of spherical sector figures, regarding their side lengths (more precisely, whether they are less, equal to greater than a quadrant) and their angles (acute, right or obtuse) and the consequences of this information on
the geometry of these figures.
\item (Chapter 4) How to find the elements (sides and angles) of a triangle from some
other given elements of this triangle
\item (Chapter 5) The Sine Rule with a large number of different proofs,
due to several mathematicians.
\item  (Chapter 6) Other trigonometric formulae, in particular formulae involving
the cosine and the tangent functions.
\item (Chapter 7) Given the sizes of three elements of a triangle
(sides or angles), how to obtain the others using trigonometry and the polar
triangle.
\end{itemize}

\medskip
\noindent{\bf An excerpt from Chapter 2.}
We now reproduce, from Chapter 2 of Book V, a discussion concerning the union of three great circles which is more extensive than the one we had in \S \ref{s:three-circles-spherical}. We reproduce it in full to give an idea of Na\d{s}\=\i r al-D\=\i n's style and his attention to details. The title of this chapter is \emph{Properties of the triangles formed by the intersections of the great circles on the surface of the sphere and their different species}. The figure is the same, and we reproduce it here (Figure \ref{fig:Trois-cercles1}). Na\d{s}\=\i r al-D\=\i n writes:

 \begin{figure}[!ht]
\centering
 \includegraphics[width=0.5\linewidth]{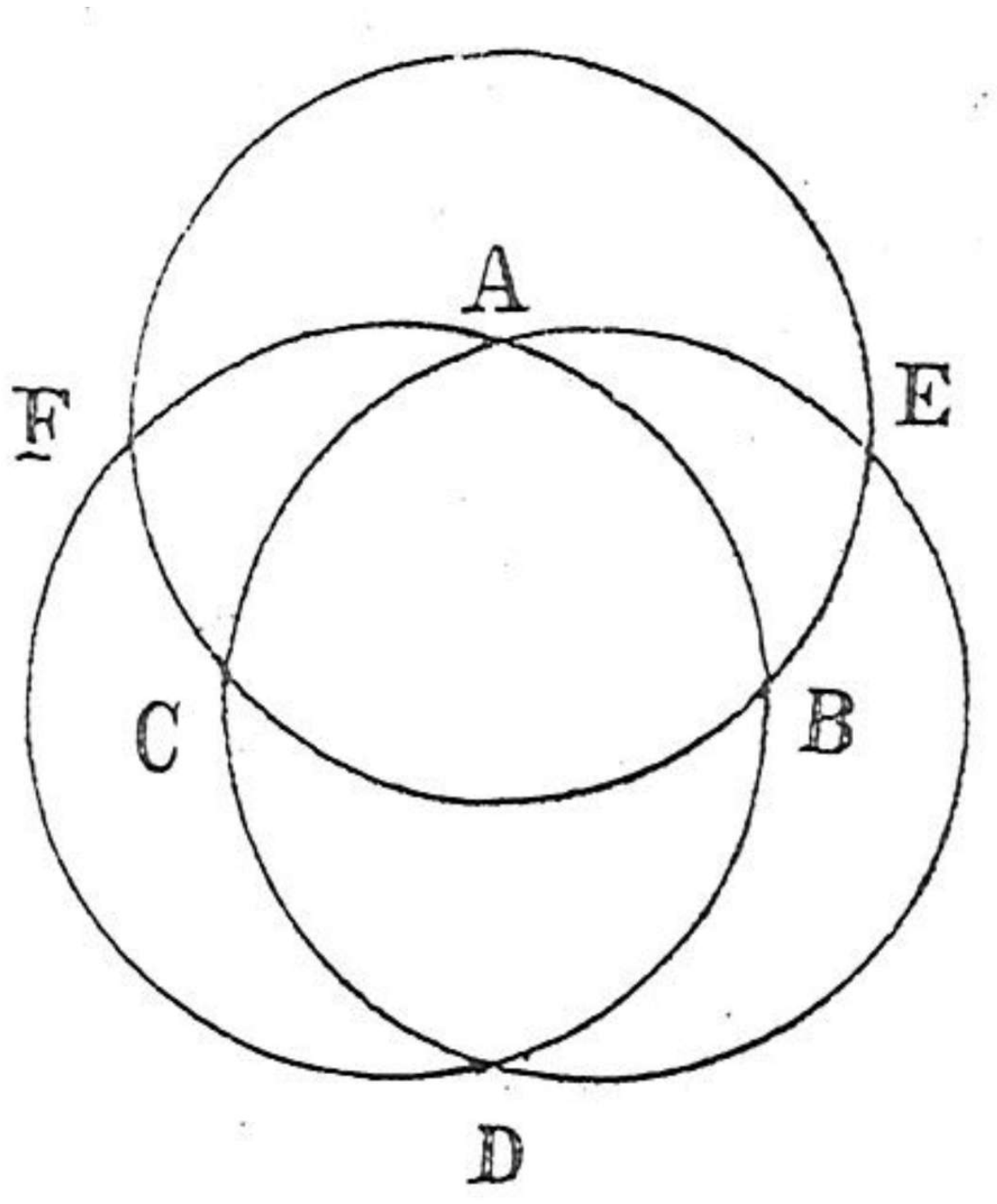}    \caption{\small 
The eight triangles made by three great circles on the sphere.}   
\label{fig:Trois-cercles1}  
\end{figure}

When three great circles intersect on the surface of the sphere, in such a way that they form a triangle, they form at the same time seven other spherical triangles; in this way the surface of the sphere is divided into eight triangles, with six intersection points, twelve arcs of circle and twenty-four angles, and as we already said, among these eight triangles, the four situated in a hemisphere are equal and similar to the four others situated in the other hemisphere, each to each.
 
 For example, the  triangle $ACF$ is equal to the triangle $DEB$  because 
 
 $AF=\frac{1}{2}\mathrm{circ.}-AB=DB;$
 
   $AC=\frac{1}{2}\mathrm{circ.}-CD=ED;$
   
     $FC=\frac{1}{2}\mathrm{circ.}-BC=EB;$
     
     $\widehat{FAC}=\widehat{EAB}$ (as they are vertically opposite, and because $\widehat{EAB}=\widehat{ABC}$) $= \widehat{EDB}$;
     
     $\widehat{AFC}=\widehat{ABC}=\widehat{EBD}$; 
     
      $\widehat{ACF}=\widehat{BCD}=\widehat{BED}$.
      
      The same will hold for the other triangles, in such a way that the four triangles situated in opposite hemispheres will be pairwise equal, for what regards the angles as well as the sides.

 For what concerns the four triangles situated in the same hemisphere, comparing any arbitrary two among them, we find that one of their  angles and two of their sides are equal, and their other elements are supplementary of each other. Thus, if we compare the triangle $ABC$ with the triangle $CDF$, we see that:
 
These two triangles have the angles at $C$ equal, as vertically equal, and the side $AB=FD$ (both being the supplement of the arc $AF$), whereas $BC=
\mathrm{suppl.} \ \widehat{FC}$; $AC= \mathrm{suppl.} \ \widehat{CD}$; the angle $\widehat{ABC}= \widehat{AFC}= \mathrm{suppl.} \ \widehat{CFD}$; 	and the angle 
$\widehat{BAC}= \widehat{BDC}= \mathrm{suppl.} \ \widehat{CDF}$.

 With this, it suffices to know only one of these eight triangles for all of them to be determined.
 
 Moreover, according to whether each side of a triangle is equal to, greater or smaller than a semicircle, we have the following ten species of triangles:\footnote{\label{f:q} In the following table and later in the text, the letter D designates a quadrant, that is, a quarter of a great circle. In his translation, Carath\'eodory writes sometimes q and other times Q. We have opted for the capital Q.}

$$A \ \  \begin{cases}
1^{\mathrm{st}} \ \mathrm{case}:  \ 3 \ \mathrm{sides} =  3\  \mathrm{Q}; \\
    2^{\mathrm{nd}}\  \mathrm{case}:  \ 2 \ \mathrm{sides}  =    \mathrm{Q};   \ 1 \ \mathrm{side}  <   \mathrm{Q};   \\
 3^{\mathrm{rd}}\  \mathrm{case}:  \ 2 \ \mathrm{sides}  =   \mathrm{Q};   \ 1 \ \mathrm{side}  >  \mathrm{Q}; 
\\   4^{\mathrm{th}}\  \mathrm{case}:  \ 1 \ \mathrm{side} =   \mathrm{Q};   \ 2 \ \mathrm{sides} <   \mathrm{Q};  \\
  5^{\mathrm{th}}\  \mathrm{case}:  \ 1 \ \mathrm{side}  = \ \mathrm{Q};   \ 2 \ \mathrm{sides}  >  \mathrm{Q}; 
\\   6^{\mathrm{th}}\  \mathrm{case}:  \ 1 \ \mathrm{side} = \  \mathrm{Q};   \ 1 \ \mathrm{side} <   \mathrm{Q};   \ 1 \ \mathrm{side} >  \mathrm{Q}; \\
 7^{\mathrm{th}}\  \mathrm{case}:  \ 3 \ \mathrm{sides}  <   \mathrm{Q};   
\\   8^{\mathrm{th}}\  \mathrm{case}:  \ 2 \ \mathrm{sides} >  \mathrm{Q};   \ 1 \ \mathrm{side} <   \mathrm{Q};  \\
  
 9^{\mathrm{th}}\  \mathrm{case}:  \ 2 \ \mathrm{sides}  <   \mathrm{Q};    \ 1 \ \mathrm{side} >  \mathrm{Q}; 
  \\   10^{\mathrm{th}}\  \mathrm{case}:  \ 3 \ \mathrm{sides} >  \mathrm{Q}.  \\
\end{cases}$$
 
 However since we find two of these species of triangles each time there is an intersection of three circles, it follows that there is a total of five species of intersections. Indeed, let us suppose that one of the eight triangles belongs to the $7^{\mathrm{th}}$ category, in other words, that each of its three sides is smaller than a quadrant. In this case, the three triangles situated in the same hemisphere belong to the $8^{\mathrm{th}}$ category since they must have necessarily two sides greater than a quadrant and one side smaller than a quadrant. Indeed, since the triangle in question should have its three sides equal to one side of each triangle, these will have a side smaller than a quadrant, whereas the two other sides will be the supplements of those of the first. Thus we can conclude from here that if the first triangle belonged to the $8^{\mathrm{th}}$ category, and if it had two sides greater than a quadrant and one side smaller, two of the remaining sides would be of that same category, whereas the third would belong to the $7^{\mathrm{th}}$ category and would have its three sides smaller than a quadrant. This is how two of the ten species of triangles that we have just enumerated (the $7^{\mathrm{th}}$ and the $8^{\mathrm{th}}$) coexist and derive from the same nature of intersection. The same can be said of the $4^{\mathrm{th}}$, the $5^{\mathrm{th}}$ and the $6^{\mathrm{th}}$ species, since among the four triangles, one will be of the $4^{\mathrm{th}}$ species, a second one of the $5^{\mathrm{th}}$, and the two others of the $6^{\mathrm{th}}$. Likewise for the $9^{\mathrm{th}}$ and for the $10^{\mathrm{th}}$ species (two of the four triangles belonging to each of these two species). Only the first species is not associated with any other species than itself.
 
 Here are now the five intersections that we have just talked about:

 $$B \ \  \begin{cases}
 
 \mathrm{The} \ 1^{\mathrm{st}} \ \mathrm{intersection \ comprises \ triangles \ of \ the}\   1^{\mathrm{st}} \  \mathrm{species}. \\

 \mathrm{The} \  2^{\mathrm{nd}}\  \mathrm{intersection \ comprises \ triangles \ of \ the} \ 2^{\mathrm{nd}}  \   \mathrm{and} \ 3^{\mathrm{rd}} \ \mathrm{species}. \\

 \mathrm{The} \  3^{\mathrm{rd}} \ \mathrm{intersection \ comprises \ triangles \ of \ the} \ 4^{\mathrm{th}}, 5^{\mathrm{th}}   \   \mathrm{and} \ 6^{\mathrm{th}}\  \mathrm{species}. \\

 \mathrm{The} \  4^{\mathrm{th}}\  \mathrm{intersection\ comprises \ triangles \ of \ the} \ 7^{\mathrm{th}}  \   \mathrm{and} \ 8^{\mathrm{th}} \ \mathrm{species}. \\

 \mathrm{The} \  5^{\mathrm{th}}\  \mathrm{intersection \ comprises \ triangles \ of \ the} \ 9^{\mathrm{th}}  \   \mathrm{and} \ 10^{\mathrm{th}} \ \mathrm{species}. \\

\end{cases}$$
 
We can also arrange the spherical triangles according to their angles, depending on whether they are right, acute or obtuse, in ten  categories as follows:\footnote{\label{f:D} In the following table, we use the letter D to mean a right angle, following Carath\'eodory who uses D for ``Droit".}

 $$A' \ \  \begin{cases}
 
1^{\mathrm{st}} \ \mathrm{case}: \  \mathrm{Each \ of \ the\  three\  angles =  D.}\\

2^{\mathrm{nd}} \ \mathrm{case}: \  \mathrm{2 \ angles =  D, 1 \ angle  <  D.} \\

3^{\mathrm{rd}} \ \mathrm{case}: \  \mathrm{2 \ angles =  D, 1 \ angle  >  D.} \\

4^{\mathrm{th}} \ \mathrm{case}: \  \mathrm{1 \ angle =  D, \ 2 \ angles  <  D.} \\

5^{\mathrm{th}} \ \mathrm{case}: \  \mathrm{1 \ angle =  D, \ 2 \ angles  >  D.} \\
  
    6^{\mathrm{th}} \ \mathrm{case}:  \mathrm{1 \ angle =D; \ 1\ angle >D; \ 1 \  angle < D.}\\

  7^{\mathrm{th}} \ \mathrm{case}:  \mathrm{ 3 \ angles  < D.}\\

  8^{\mathrm{th}} \ \mathrm{case}:  \mathrm{ 1 \ angle  < D, \  2 \ angles  > D.}\\

  9^{\mathrm{th}} \ \mathrm{case}:  \mathrm{ 3 \ angles  > D.}\\

 10^{\mathrm{th}} \ \mathrm{case}:   \  \mathrm{1 \ angle >  D, \ 2 \ angles  <  D.} 
\\

\end{cases}$$

The number of intersections giving rise to these last ten categories is also five:

 $$B' 
 \ \  \begin{cases}
 \mathrm{The} \ 1^{\mathrm{st}} \ \mathrm{intersection \ comprises  \ triangles \  of\  the} \ 1^{\mathrm{st}}\  \mathrm{\ species} \\
\mathrm{The} \ 2^{\mathrm{nd}} \ \mathrm{intersection \ comprises  \ triangles \  of\  the} \ 2^{\mathrm{nd}} \ \mathrm{and} \  3^{\mathrm{rd}} \ \mathrm{species} \\

     \mathrm{The} \ 3^{\mathrm{rd}} \ \mathrm{intersection \ comprises  \ triangles \  of\  the} \ 4^{\mathrm{th}}, \ 5^{\mathrm{th}}  \ \mathrm{and} \  6^{\mathrm{th}} \ \mathrm{species} \\
     
        \mathrm{The} \ 4^{\mathrm{th}} \ \mathrm{intersection \ comprises  \ triangles \  of\  the} \ 7^{\mathrm{th}},  \mathrm{and} \  8^{\mathrm{th}} \ \mathrm{species} \\

        \mathrm{The} \ 5^{\mathrm{th}} \ \mathrm{intersection \ comprises  \ triangles \  of\  the} \ 9^{\mathrm{th}} \ \mathrm{and} \  10^{\mathrm{th}} \ \mathrm{species} \\

\end{cases}$$
 
For what regards the way of finding the association of the various categories of triangles noted under each intersection species,  we shall adopt an approach which is similar in every point to the one we used while we treated the first five intersections (B) relative to the ten categories of triangle $(A)$, ordered according to the magnitude of their sides.
 
\medskip
\noindent{\bf An excerpt from Chapter 3.}
 Without going into the proofs, we quote now a few propositions from Chapter 3 of Book V. These propositions are more involved than those of Chapter 2. This chapter is titled \emph{On the rules relative to the various species of triangles and on the general and particular considerations attached to them}.

\medskip
\noindent {\bf Proposition I.} The angles of any triangle whose sides are quadrants are necessarily equal to 90$^{\mathrm{o}}$. The intersections of the sides constitute [coincide with] the poles of the sides, each vertex is the pole the opposite side. 

\medskip
\noindent {\bf Proposition II.}   If a spherical triangle has two sides equal to a quadrant, and the third side  less than a quadrant, then it will have  two angles equal to  90$^{\mathrm{o}}$ and an acute angle. The vertex of the latter is the pole of the arc that is opposite to it, and the poles of the other two sides are placed on the arc of the acute angle, outside the triangle.

\medskip
\noindent {\bf Proposition III .}  If two sides of a spherical triangle are equal to a quadrant and the third greater than a quadrant, then this triangle has two right angles and an obtuse angle (the one opposite to the side greater than a quadrant). The vertex of the latter is the pole of the arc that is opposite to it. As for the two other sides, their poles lie on the arc of the obtuse angle in the interior of the triangle.

\medskip
\noindent {\bf Proposition IV.}  If a spherical triangle has one side equal to a quadrant and the two others smaller than a quadrant, then the angle opposite to the quadrant is obtuse, and the two other angles are acute. Furthermore, the three poles will lie inside the triangle.

\medskip
\noindent {\bf Proposition V.}  If a spherical triangle has one side equal to a quadrant and the two others  greater than a quadrant, then all its angles are obtuse. The poles will be situated in the interior of the triangle.

\medskip
\noindent {\bf Proposition VI.}  If a spherical triangle has one side equal to a quadrant, one side smaller and one side greater than a quadrant, then the angle opposite to the greatest side is obtuse and the two other angles are acute. The three poles fall outside the triangle.

\medskip
\noindent {\bf Proposition VII.}  If the three sides of a triangle are less than a quadrant, then the triangle has two acute angles. The third angle can be right, obtuse or acute. The three poles of this triangle fall outside the triangle.

\medskip
\noindent {\bf Proposition VIII.}  If a triangle has two sides greater  than a quadrant and the third one smaller than a quadrant, then its angles are necessarily of one of the following types:
\begin{enumerate}

\item One right angle and two obtuse angles.

\item One right angle, one acute angle and one obtuse angle.

\item One acute angle and two obtuse angles.

\item One obtuse angle and two acute angles.

\item Three obtuse angles.

\end{enumerate}

\medskip
\noindent {\bf Proposition IX.} If a triangle has one side greater than a quadrant and the two others smaller than a quadrant, then the angle opposite to the greatest side is obtuse and the two others are acute. The three poles fall outside the triangle.

\medskip
\noindent {\bf Proposition X.}  If  the three sides of a triangle are all greater than a quadrant, then the three angles are obtuse and the three poles fall in the interior of the triangle.

Chapters  4 to 6 of Na\d{s}\=\i r al-D\=\i n's treatise are dedicated to the proofs of the spherical trigonometric formulae. The proofs use the notion of polar triangle.  

\medskip
\noindent{\bf An excerpt from Chapter 5.}
 Chapter 5 is titled ``On the so-called supplementary figure and its different species".
 Na\d{s}\=\i r al-D\=\i n mentions in this chapter a large number of mathematicians who contributed to the formulae of spherical trigonometry. They include in particular 
   Abu al-Ray\d{h}\=an al-B\=\i r\=un\=\i,\index{Biruni@al-B\=\i r\=un\=\i}
  Ab\=u Na\d{s}r Man\d{s}\=ur Ibn `Ir\=aq,\index{Ibn Iraq@Ibn `Ir\=aq}
   Ab\=u al-Waf\=a'  Mu\d{h}ammad Ibn-Mu\d{h}ammad al-B\=uzj\=an\=\i,\index{Buzjani@al-B\=uzj\=an\=\i}
   Ab\=u Ma\d{h}m\=ud \d{H}am\=\i\ d ibn al Azhar  al-Khujand\=\i,\index{Khujandi@al-Khujand\=\i}
  Abul Fa\d{d}l al-Nayr\=\i z\=\i,\index{Nayrizi@Al-Nayr\=\i z\=\i}
  in his commentary on the \emph{Almagest},
   Ab\=u Ja`far al Kh\=azin\index{Khazin@al Kh\=azin}
    in his book titled \emph{Partial researches on the inclination of the inclinations and introduction of the straight sphere},
   K\=ushy\=ar ibn Labb\=an al-J\=\i l\=\i,\index{Jili@al-J\=\i l\=\i}
   Th\=abit ibn Qurra\index{Thabit@Th\=abit ibn Qurra} and several others.
     The history of spherical trigonometry in the medieval Arab world has been investigated relatively recently, and all the authors mentioned by  Na\d{s}\=\i r al-D\=\i n play some role in this history, see \cite{Bellosta, Debarnot1, Debarnot2, R1}.

    The chapter starts with the following central principle, namely, the sine rule:
 
 \medskip
 
\noindent {\bf Proposition.} The ratios of the sines of the sides of the triangles formed on the surface of a sphere by the intersection of arcs of great circles are equal to the ratios of the sines of the angles opposite to these sides.
 
 \medskip

  Na\d{s}\=\i r al-D\=\i n writes: ``We usually start by establishing this principle for what concerns the right triangle and to this effect we follow various ways that are all exposed in the book of the scholar Abu al-Ray\d{h}\=an al-B\=\i r\=un\=\i \ titled  \emph{The keys of the knowledge of the figures of the  spherical and other surfaces}".\footnote{See Al-B\=\i r\=un\=\i, \emph{Kit\=ab Maq\=al\=\i d `Ilm al-hay'a (The Keys to the Science of Astronomy): La trigonom\'etrie sph\'erique chez les Arabes de l'Est \`a la fin du Xe si\`ecle},   Marie-Th\'er\`ese Debarnot (editor and translator), Damascus, 1985.}  

He then gives    proofs of this special case, each time mentioning their author (Ab\=u Na\d{s}r,  Ab\=u al-Waf\=a', al-Khujand\=\i  , Abu al-Ray\d{h}\=an,  etc.)  One of these proofs uses the sector figure. 
  He deduces from this several trigonometric formulae.

  At the end of Chapter 5,  Na\d{s}\=\i r al-D\=\i n  writes the following:
 
    \begin{quote}\small
    Ab\=u Ma\d{h}m\=ud al-Khujand\=\i\ gave to this figure [the supplementary figure] the name ``astronomy rule". Others called it \emph{The figure which exempts from the quadrilateral}.  In his book titled \emph{The keys of the knowledge of what happens on the surface of the sphere},  Abu al-Ray\d{h}\=an declares that it is indeed the Emir Ab\=u Na\d{s}r who was the first to use the supplementary figure instead of the quadrilateral, but that the name which it carries comes from   K\=ushy\=ar ibn Labb\=an al-J\=\i l\=\i . However, this assertion comes up against a difficulty, because the Emir  Ab\=u Na\d{s}r, in the second part of the first book of the work titled \emph{The Royal Almagest}  at the beginning of Chapter III, where this figure is mentioned, writes verbatim: ``The third chapter of what can exempt from the quadrilateral figure"; and then, after he has mentioned the treatise of Th\=abit ibn Qurra on the different varieties that occur in the quadrilateral figure, he adds: ``And  Th\=abit ibn Qurra also wrote a treatise on what exempts from the quadrilateral figure, but he who utilizes it must know the use of the\index{compounded ratio} compounded ratios; but I will show right here a procedure which exempts you from the quadrilateral as well as from the compounded ratios." These words prove indeed that the term ``supplementary figure" itself is due to the Emir Ab\=u Na\d{s}r  who received it from Th\=abit ibn Qurra.\footnote{Carath\'eodory writes in a footnote that a note in the margin of the manuscript reads: ``I add: That the Emir  Ab\=u Na\d{s}r talked about this in his \emph{Commentary on Menelaus}, as I have already reported on, on the occasion of the Supplementary Figure."}
    
    \end{quote}
    
        We end this paper by quoting a few results from Chapter 6.
        
\medskip
\noindent{\bf A few results from Chapter 6.}
  This chapter is titled ``On the figure called shadowed, on its consequences and its accessories".
  
  The shadow is what we call the tangent of an angle. Again, in addition to the mathematics, Na\d{s}\=\i r al-D\=\i n provides us with historical information. He writes that 
the priority, for what concerns this figure, is due to Ab\=u al-Waf\=a'  al-B\=uzj\=an\=\i , an information provided by Abu al-Ray\d{h}\=an. 
     He states the following proposition which concerns right spherical triangles:
     
          \[\frac{\mathrm{the \ sine\ of \ one\ of \ the \ sides \ of\ the \ right\ angle}}{\mathrm{the\ sine \ of\ the\ right \ angle}}\]
     \[=
     \frac{\mathrm{the \ shadow \ of \ the \ other \ side \ of\ the \ right\ angle}}
   {\mathrm{the \ shadow \ of\ the  \ angle \ opposite \ to \ this \ side}}.\]

From this proposition, Na\d{s}\=\i r al-D\=\i n deduces several corollaries, among which we state the following:

\medskip
\noindent {\bf Corollary.}   In a right triangle, the cosine of the acute angle that we assume it has is to the the sine of the given right angle like the tangent of the complement of the arc opposite to this right angle to the tangent of the complement of the side contained between the right angle and the angle  that we assumed to be acute, or else as the tangent of the side contained by the two angles is to the tangent of the side opposite to the given right angle.

\medskip

 Na\d{s}\=\i r al-D\=\i n  ends this chapter by several remarks on the tangent function,  saying that one shortage is the rapid increase of the tangent of the arcs that exceed 45$^{\mathrm{o}}$. He writes: ``The tangents increase rapidly beyond the radius, because already the tangent of the eighth part of the circle (45$^{\mathrm{o}}$) is equal to the radius. Thus if you write the tangents in a table in which the arcs increase by equal degrees, the differences between the tangents beyond 45$^{\mathrm{o}}$ become very large.
This is why we do not have much confidence in taking these tables and modifying what is actually between the lines, as it is done for the other tables.
However, if we divide the arcs (in two categories), we will find that the objection that we made on this account to this figure is groundless, since it is not necessary to take tangents in the tables. 
In any circumstance, we will be able to operate by using the tangents of the eighth part of the circle, because the tangents have properties that we do not find in the sines, and which allow us to use some of them for the others. We have said a few words about this matter at the beginning of each chapter, and we shall soon show ow we can use all the tangents by knowing only the tangents of the arcs that are less than the eighth part of the circle."

    \bigskip
 We have reproduced these excerpts from Book V of 
 Na\d{s}\=\i r al-D\=\i n al-\d{T}\=us\=\i 's  treatise, written before the introduction of the trigonometric functions in Western Europe, to give an idea of its content, both from the mathematical and from the historical point of view.
We have also tried to suggest the preciseness of Na\d{s}\=\i r al-D\=\i n's writing and his aim to be exhaustive:  he always considers all the cases of a construction or a proposition.

 \bigskip
 
Author's address: 
Athanase Papadopoulos,  Institut de Recherche Math\'{e}matique Avanc\'{e}e,   
 Universit\'{e} de Strasbourg et CNRS,   
  7 rue Ren\'{e} Descartes,  67084 Strasbourg Cedex, France.
  email: athanase.papadopoulos@math.unistra.fr


\begin{thebibliography}{99}
   
   \bibitem{Bellosta} H. Bellosta, Le trait\'e de Th\=abit ibn Qurra sur la figure secteur, 
Arabic Sciences and Philosophy 14 (2004) No. 1, p.145--168 
   
   \bibitem{Carra} B. Carra de Vaux,  Les sph\`eres c\'elestes selon Na\d{s}\=\i r al-D\=\i n al-\d{T}\=us\=\i. In: P. Tannery, Recherches sur l'histoire de l'astronomie ancienne, Gauthier-Villars, Paris, 1893, Appendice VI, p. 337-61.

\bibitem{Debarnot1} M. T. Debarnot, La trigonom\'etrie sph\'erique chez les Arabes de l'Est \`a la fin
 du Xe si\`ecle, Damascus, 1985.
 
\bibitem{Debarnot2} M. T. Debarnot, Trigonometry, In: Encyclopedia of the history of arabic science, 3 volumes, ed. R. Rashed, Routledge, London and New York, 1996, 
Vol. 2, p. 470-513.




\bibitem{Halma} N. Halma (ed.),  Composition math\'ematique de Claude Ptol\'em\'ee, traduite par Nicolas Halma et suivie
des notes de Jean Baptiste Delambre, 2 vols., Paris, Henri Grand, 1813.
Reprint, 1988, Librairie Blanchard, Paris.

    \bibitem{Euclid-Heath} T. L. Heath, The thirteen books of Euclid's Elements, 3 volumes, Cambridge University Press, 1908, Reprint, Dover.
     
   \bibitem{Heath-History3} T. L. Heath, A history of Greek mathematics, 2 volumes, Clarendon Press,  Oxford, 1921, Reprint, Dover.
 \bibitem{K1} E. S. Kennedy, Two Persian astronomical treatises by Na\d{s}\=\i r al-D\=\i n al-\d{T}\=us\=\i,  Centaurus 27 (2) (1984), p. 109--120.

 \bibitem{K2} A. Kubesov, The commentaries of Na\d{s}\=\i r al-D\=\i n al-\d{T}\=us\=\i \ on the treatise of Archimedes \emph{On the sphere and cylinder} (Russian), Voprosy Istor. Estestvoznan. i Tehn. Vyp. 2 (1969), p. 23--28.

           

 \bibitem{Papa-Theret-Lambert} A. Papadopoulos and G. Th\'eret, Johann Heinrich Lambert's memoir  \emph{Theorie der Parallellinien}: a review with commentary, 
 this volume, Chapter 7.
 
     

      \bibitem{RR2} R. Rashed and A. Papadopoulos, Menelaus'  Spherics:  Early Translation and al-M\=ah\=an\=\i /al-Haraw\={\i}'s Version,
(Critical edition of Menelaus' Spherics from the Arabic manuscripts, with historical and mathematical commentaries), De Gruyter, Series: Scientia Graeco-Arabica,  21,  2017, 890 pages.



\bibitem{RR} R. Rashed (ed.) Encyclopedia of the history of Arabic science, 3 volumes, Routledge, 1996.



   \bibitem{R1} B. A. Rosenfeld, On the mathematical works of Na\d{s}\=\i r al-D\=\i n al-\d{T}\=us\=\i \ (Russian), Istor.-Mat. Issled. 4 (1951), p. 489--512.
   
   \bibitem{Rosenfeld-Yoush} B. A. Rosenfeld et A. P. Youshkevich,  Geometry.  In  Encyclopedia of the history of Arabic science,  ed. R. Rashed et R. Morelon, Vol. 2,  Routledge, London, 1996, p. 447-494.



\bibitem{Ragep} F. J. Ragep, Na\d{s}\=\i r al-D\=\i n al-\d{T}\=us\=\i's Memoir on Astronomy, 2 volumes, Springer-Verlag New York, 1993.
   \bibitem{Rosenfeld} B. A. Rosenfeld,    A History of Non-Euclidean Geometry: Evolution of the Concept of a Geometric Space. Translated by A. Shenitzer,
   Springer-Verlag New York, 1988.  First Russian edition, 1976. 

  \bibitem{RY} B. A. Rosenfeld and A. P. Yushkevich, Notes on the treatise of  Na\d{s}\=\i r al-D\=\i n al-\d{T}\=us\=\i\ on parallel lines (Russian), Istor.-Mat. Issled. 13 (1960), p. 525--532.
  
   \bibitem{Strauss} Johann Strauss, 
Le livre fran\c cais d'Istanbul (1730-1908), Revue des mondes musulmans et de la M\'editerrann\'ee, 
 87-88 (1999), p. 277-301



    \bibitem{Tusi}  Trait\'e du
quadrilat\`ere, attribu\'e \`a Nassiruddin-el-Toussy, d'apr\`es un manuscrit tir\'e de la
biblioth\`eque de S. E. Edhem Pacha, ancien grand vizir, traduit par Alexandre Pacha
Carath\'eodory, ancien ministre des affaires \'etrang\`eres, Constantinople, Typographie et
Lithographie Osmani\'e, 1891
 

 \bibitem{Tusi-Hyderabad} Na\d{s}\=\i r al-D\=\i n al-\d{T}\=us\=\i , Majmu`al-ras\=a'il, 2 vol., 
 Hyderabad, Osmania
Oriental Publication Bureau.

 
 \end{thebibliography}
\end{document}